\documentclass{article}
\usepackage{amsmath}
\usepackage[russian]{babel}

\def\R{{\mathbf R}}
\def\Z{{\mathbf Z}}
\def\div{{\rm \,div\,}}

\newcommand{\myref}[1]{(\ref{#1})}

\def\mysection#1{\section{#1}\setcounter{equation}{0}}

\newtheorem{theorem}{Theorem}

\newtheorem{proposition}[theorem]{Proposition}

\title{Irregular dynamic systems according to R.J. DiPerna and P.L. Lions}
\author{Alexander Ovseevich}
\date{}
\begin{document}\large
\maketitle 
\noindent{\bf Abstract} These are  notes of a seminar held at the Institute for Problems in
Mechanics, RAS in 2003 and aimed at presentation of \cite{dl}. We discuss the notion of a
generalized solution to a singular Ordinary Differential Equation introduced by DiPerna and Lions.
We stress importance of singular dynamic systems from the ``philosophia naturalis'' point of view,
and extend and simplify the original approach by R.J. DiPerna and P.L. Lions. Further extensions are
discussed.

\medskip


\medskip


\mysection{Singular dynamic systems and {``philosophia na\-tu\-ra\-lis''}} In Newtonian mechanics
the Universe is governed by
\begin{equation}
m\ddot{x}=\frac{\partial U(x)}{\partial x}  \label{1}
\end{equation}
where $U$ is the potential energy. Newton himself discovered that the gravitation cor\-res\-ponds to
\begin{equation}
U(x)=-G\sum_{i\neq j}\frac{m_{i}m_{j}}{|x_{i}-x_{j}|},  \label{2}
\end{equation}
where $G$ is a  constant, $m_{i}$ is the mass of the $i$th particle, and $|y|$ stands for the length
of 3-dimensional vector $y$.

If one takes this this ODE point of view seriously,  mathematical facts are to be regarded as
philosophical principles. E.g. the Laplace determinism, being a belief that the present determines
future, is modelled by a uniqueness theorem for the Cauchy problem for (\ref{1}). Similarly, the
existence theorem also has a physical or philosophical meaning as a claim that any present state has
a future (however unfavorable it can be). So, the existence and uniqueness theorem for ODE might
be of more than a purely mathematical interest. Unfortunately, the
standard existence and uniq\-ue\-ness theorem for ODE 
(under the Lipschitz con\-di\-tion for the force function $\frac{\partial U(x)}{%
\partial x}$)
is not applicable, say, to potential energy ($\ref{2}$), for the corresponding forces are not
Lipschitz, and not even continuous.

Physical arguments due to Boltzmann and Loschmidt make the issue of existence and
uniqueness for (\ref{1}) rather problematic. For starters, the Newton law (%
\ref{1}) is (or at least seems to be) reversible in time: if $t\mapsto x(t)$ is a solution, then
$t\mapsto x(-t)$ is a solution as well. However, the Universe has entropy which increases in time
according to the second law of thermodynamics.   It is clear, that if the phase flow is reversible,
the entropy should oscillate, and cannot be monotone increasing.

One can imagine, however, that the reversibility does not hold: the ``solution'' $x(-t)$ is not, in
fact, a solution. Indeed, the trivial ``proof'' that it is a solution is a formal application of
rules of differentiation. However,  the ``true'' solution $t\mapsto x(t)$ might obey (\ref{1}) in a
generalized sense, and  may not be differentiable. Then, it is possible that $x(-t)$ does not
satisfy (\ref{1}) in the generalized sense.

These remarks suggest  that the search for a proper notion of a solution of a singular ODE and a
further study of its properties rather belongs to natural philosophy than to a routine mathematics
and is totally justified.

In fact, this issue is relevant well beyond reconciliation of classical mechanics and
thermodynamics. For instance, a natural source of singular dynamic systems is the Control Theory,
where the Pontryagin Hamiltonians usually are not everywhere dif\-fe\-ren\-ti\-able, and the
corresponding vector fields have jumps.

\mysection{Cauchy-Lipschitz theorem} This is a well-known theorem on existence and uniqueness of
solutions to the Cauchy problem for ODE. It is so well-known that mathematicians usually do not
suspect that anybody need something else in this area. Its statement is as follows:
\begin{theorem}\hglue -.3em{\bf .}
Suppose we are given ODE
\begin{equation}\label{ode}
  \dot x=b(x), \, x\in\mathbf{R}^n,\, b\in\mathbf{R}^n,
\end{equation}
where $b$ is a Lipschitz continuous vector field. Then, for any $x_0\in\mathbf{R}^n$ there is a
unique $C^1$ solution to \myref{ode} such that $x(0)=x_0$. Moreover, if the phase flow
$X_t:\mathbf{R}^n\to\mathbf{R}^n$ is defined by $X_t(x_0)=x(t)$ then the phase flow is Lipschitz
continuous w.r.t. the variable $x_0\in\mathbf{R}^n$, and  $C^1$ w.r.t. the variable
$t\in\mathbf{R}$.
\end{theorem}
Lipschitz functions are exactly the functions of the Sobolev class $W^{1,\infty}$ with uniformly
bounded generalized derivative. In other words, if $b$ belongs to $W^{1,\infty}$, then the
cor\-res\-pon\-ding ODE has good pro\-per\-ties. In DiPerna--Lions paper \cite{dl} it is shown that
rich theory exists for $b\in W^{1,1}$. Recall, that a function $u$ which locally belongs to $L^1$ is
said to belong to $W^{1,p}$ if $u\in L^p$, and its first derivative $\frac{\partial u}{\partial x}$
in the sense of distributions belongs to $L^p$. Here, $p\in[1,\infty]$. We also use the space
$W_*^{1,1}$ of $L^1$ functions such that their distributional first partial derivatives are bounded
measures. This latter space can be regarded as a ``weak version'' of $W^{1,1}$.

In this paper we extend the DiPerna--Lions theory to  vector fields $b\in W_*^{1,1}$. This is
important because most commonly used singular vector fields like that with jump singularities along
a hypersurface belong to this class, and do not belong to the original DiPerna--Lions class
$W^{1,1}$. \mysection{Extended DiPerna--Lions theory} We build our  exposition for simplicity not
around differential equa\-t\-ions in $\mathbf{R}^n$, but around differential equations on a torus
${\mathcal T}=\mathbf{R}^n/\mathbf{Z}^n$. Any other closed manifold is as good as torus for our
purposes, but in the torus case we can utilize almost the same classical notations as in the
euclidean case. First, state the {\it extended DiPerna--Lions conditions} on the vector field $b$:
\begin{equation}\label{dip_l}
    \div b\in L^\infty({\mathcal T}), \quad
    b\in W_*^{1,1}({\mathcal T}).
\end{equation}
The original DiPerna--Lions conditions were stated in the Euclidean setting and require that
\begin{equation}\label{dip_le_a}
    \div b\in L^\infty(\mathbf{R}^n), \quad
    b\in W^{1,1}_{\rm loc}(\mathbf{R}^n),\quad
\end{equation}
\begin{equation}\label{dip_le_b}
    \frac{b(x)}{1+|x|}\in L^\infty(\mathbf{R}^n)+L^1(\mathbf{R}^n)
\end{equation}
which is clearly more involved. We can extend the DiPerna--Lions theory to the Euclidean setup by
requiring
\begin{equation}\label{dip_le2a}
    \div b\in L^\infty(\mathbf{R}^n), \quad
    b\in {W_{*{\rm loc}}^{1,1}}(\mathbf{R}^n),\quad
\end{equation}
\begin{equation}\label{dip_le2b}
    \frac{b(x)}{1+|x|}\in L^\infty(\mathbf{R}^n)+L^1(\mathbf{R}^n)
\end{equation}
instead of \myref{dip_le_a}, \myref{dip_le_b}.

 The only
difference is the replacement of the ``weak'' space $W_*^{1,1}$ with its ``strong'' version
$W^{1,1}$.

A typical example of a singular $b$ which satisfies the DiPerna--Lions conditions \myref{dip_le_a}
is the Hamiltonian field with the Hamiltonian function
\begin{equation}\label{ham}
  H(p,q)=\sum_i\frac{p_i^2}{2m_i}+\frac{1}{2}\sum_{i\neq
  j}\frac{e_{i}e_{j}}{|q_{i}-q_{j}|^\alpha},
\end{equation}
where $q_i\in{\mathbf R}^3$ and $\alpha\in(0,1)$. Note that the Coulomb system ($\alpha=1$) does not
satisfy both \myref{dip_le_a} and \myref{dip_le2a}. The growth condition \myref{dip_le_b} does not
hold for  the Hamiltonian field \myref{ham}.

\if{false} Strictly speaking, to satisfy the last requirement $ \frac{b(x)}{1+|x|}\in L^\infty+L^1$
it is necessary to confine ourselves to the invariant subspace $$ \sum p_i=0,\quad \sum m_iq_i=0.
$$
\fi
 The ``less singular'' Hamiltonian $H(p,q)=|p|$, which is
typical for the control theory,  does not fit the original DiPerna--Lions conditions
\myref{dip_le_a}, but fits \myref{dip_le2a}.

In fact, there is no such a thing as the DiPerna--Lions theorem parallel to that of
Cauchy-Lipschitz. What does exist is the DiPerna--Lions theory, which is only partially concerned
with ODE.
\subsection{Transport equation}
This is the equation
\begin{equation}\label{teq}
  \frac{\partial u}{\partial t}=\sum_i\,b_i \frac{\partial u}{\partial
  x_i},
\end{equation}
which is dual to the conservation law equation
\begin{equation}\label{dens}
  \frac{\partial \rho}{\partial t}+\div b\rho=0
\end{equation}
 describing evolution of the density
of particles moved by the phase flow of \myref{ode}. If the vector field $b$ is sufficiently
regular, say, if the Cauchy-Lipschitz condition holds, one can write down the general solution of
\myref{teq} in terms of the phase flow. Namely,
\begin{equation}\label{sol}
  u(x,t)=u^0(X_t(x)),
\end{equation}
where the function $u^0(x)=u(x,0)$ is arbitrary. In other words, the phase flow of ODE defines and
is defined simultaneously by the solution of the Cauchy problem for the transport equation
\begin{equation}\label{cp}
\frac{\partial u}{\partial t}=b\cdot\nabla u,\quad u(x,0)=u^0(x).
\end{equation}
\subsection{Renormalizable and approximable solutions}
The approach adopted by DiPerna--Lions is to study the Cauchy problem \myref{cp} for the transport
equation without recourse to the corresponding ODE, and then {\it define} the phase flow via
\myref{sol}.

The first step is to define what is the solution to \myref{cp}. Of course we have a notion of the
classical solution: a differentiable function $u$ which satisfy \myref{cp}. However, one cannot
expect to solve \myref{cp} in the classical sense if the vector field $b$ is not sufficiently
regular. The correct definition of the solution is achieved in two steps. First, we recall the old
and well known notion of the {\it weak} solution. A function $u(x,t)$ is a weak solution of
\myref{cp} if for every smooth function $\phi$ with compact support in ${{\mathcal T}}\times [0,T )$
(test function) we have
\begin{equation}\label{ws}
\int_0^T dt \int dx\, u\frac{\partial \phi}{\partial t}=-\int dx\, u^0(x)\phi(x,0)+\int dx\,
u\div(b\phi).
\end{equation}
In other words, we multiply the formal equality \myref{cp} by $\phi$ and formally integrate by part.
One can prove the existence of the weak solution of \myref{cp} under very weak assumptions on the
vector field $b$. E.g., it suffices to assume that $b\in L^1$ and $\div b\in L^1$. One cannot,
however, guarantee the uniqueness of the weak solution. To restore the uniqueness DiPerna and Lions
invented a new notion of {\it re\-nor\-ma\-li\-z\-able} solution. This requires a new set of test
functions. Suppose that $\beta:\mathbf{R}\to\mathbf{R}$ is a $C^1$-function which is bounded itself
and has bounded derivative.

A function function $u(x,t)$ is a {\it re\-nor\-ma\-li\-z\-able} solution of \myref{cp} if for each
above $\beta$ the function $\beta(u)$ is a weak solution of \myref{cp}. It is clear that if $u$ is a
classical solution, then $\beta(u)$ also is, but for the weak solutions this trans\-for\-ma\-tion
may fail to give a solution. In fact, this notion of renormalizability is very close to well known
entropy conditions for solutions of {\it nonlinear} equations of con\-ser\-va\-tion laws
\cite{waves}, \cite{kruzh}.

However, for the extended DiPerna--Lions conditions \myref{dip_l} we find it more appropriate to
define and work with another type of solutions --- {\it approximable} solutions. We say that a
function $u\in L^\infty({\mathcal T}\times[0,T])$ is an approximable solution of \myref{cp} if   $u$
is a  limit $u_{\epsilon}(t)\to u(t)$ in $H^{-s}$ uniformly in $[0,T]$ of functions $u_\epsilon$,
which are smooth w.r.t. the space variables and satisfy
\begin{equation}\label{cp_epsilon}
\frac{\partial u_{\epsilon}}{\partial t}=\sum_i\,b_i \frac{\partial u_{\epsilon}}{\partial
  x_i}+r_{\epsilon},\quad u_{\epsilon}(x,0)=u^0_{\epsilon}(x),
\end{equation}
where $r_{\epsilon}$ are measures w.r.t. $x$ such that their total variations $\Vert
r_{\epsilon}^t\Vert$ are uniformly bounded for each $t\in[0,T]$ and tend to zero as $\epsilon\to~0$.
Note that J. Moser \cite{M} stressed that for numerous
 problems approximate solutions  can be
more valuable then the exact ones. Our main theorem (extended DiPerna--Lions theorem) is about
existence and uniqueness of approximable solution to \myref{cp}.
\begin{theorem}\hglue -.3em{\bf .}\label{dlth}
Suppose that the extended DiPerna--Lions conditions \myref{dip_l} hold, and $u^0\in
L^\infty({\mathcal T})$. Then there exist a unique approximable solution $u$ to \myref{cp}. This
solution is renormalizable and belongs to
\begin{equation}\label{space}
  u\in L^\infty(0,T;L^\infty({\mathcal T}))\cap C([0,T];L^p({\mathcal T}))
\end{equation}
for each $1\leq p <\infty$.
\end{theorem}
In other words, the solution $u$ with bounded initial condition is bounded and depends on time $t$
in a continuous way. Denote by $T_t$ the Cauchy operator
\begin{equation}\label{Cauchy}
  T_t(u^0)=u^t,
\end{equation}
where $u^t(x)=u(x,t)$. By using Theorem \ref{dlth} one can restore the phase flow. This requires a
general result from functional analysis.
\begin{theorem}\hglue -.3em{\bf .}\label{dpth}
Suppose that
\begin{equation}\label{hom}
  A: L^\infty({\mathcal T})\to L^\infty({\mathcal T})
\end{equation}
is a (automatically continuous) homomorphism of rings with unit ($A(fg)=A(f)A(g)$, $A1=1$). Then,
$A$ is a measurable change of variables: there exists $\Phi:{\mathcal T}\to{\mathcal T}$ such that
$Au(x)=u(\Phi(x))$.
\end{theorem}
One can see easily from the definition of an approximable
solution that $T_t(u^2)=T_t(u)^2$ (here, $T_t$ the Cauchy operator and $u^2=u\times u$), and,
therefore, one can apply Theorem \ref{dpth} to $T_t$. We obtain
\begin{equation}\label{flow}
T_t(u)(x)=u(X_t(x)),
\end{equation}
where $X_t$ is a one parameter group of measurable transformations of ${\mathcal T}$. This is the
phase flow we were looking for.

Notice that our construction of the phase flow requires studying the Cauchy problem \myref{cp} only
for {\it bounded} initial data. At that point the original approach of DiPerna--Lions is different.
They have built the flow $X_t:{\mathbf R}^n\to{\mathbf R}^n$ as a renormalizable solution
 to
 \begin{equation}\label{flow2}
\frac{\partial X(x,t)}{\partial t}=\sum_i\,b_i \frac{\partial X(x,t)}{\partial
  x_i},\quad X(x,0)=x,
\end{equation}
where the initial data is surely unbounded. In our approach this primary motivation for introducing
and studying renormalizable solutions disappear. Basically by the same reason we do not study the
Cauchy problem \myref{cp} with initial data from $L^p, \,p<\infty$.

Another important aspect of our extended DiPerna--Lions theory is the stability the\-o\-rem for
approximable solutions.
\begin{theorem}\hglue -.3em{\bf .}\label{stabth}
Suppose that vector fields $b_n\in L^1({\mathcal T})$ are such that $b_n, \div(b_n)$ converge in
$L^1$ to (respectively) $b, \div b$, where $b$ satisfy the extended DiPerna--Lions conditions
\myref{dip_l}. Suppose also, that $u_n$ is a boun\-ded sequence in $L^\infty(0,T;L^\infty)$ of
approximable solutions of \myref{cp} with $b$ replaced by $b_n$, and assume that $u_n^0\to u^0$ in
$L^1$. Then, $u_n$ converges as $n\to\infty$ in $C([0,T];L^1)$ to the approximable solution of
\myref{cp} corresponding to the initial condition $u^0$.
\end{theorem}
In terms of the phase flow this means that  disturbances of the vector fields which are small in
$L^1$ and produce small disturbances of the divirgences in $L^1$, give a small change of the flow.
\mysection{Open questions} We mention only a few arbitrarily chosen issues.
\subsection{Formally reversible system with irreversible dynamics}
The Coulomb system does not satisfy DiPerna--Lions conditions \myref{dip_l}, and
hypo\-the\-ti\-cally, in general, there is no phase flow and re\-nor\-ma\-li\-z\-able solutions in
the sense of DiPerna--Lions. One can expect, though, that it is possible to solve the corresponding
Cauchy problem by using the vanishing viscosity method. In other words, we are going to solve the
Cauchy problem
\begin{equation}\label{visc}
\frac{\partial u_\epsilon}{\partial t}=\sum_i\,b_i \frac{\partial u_\epsilon}{\partial
  x_i}+\epsilon^2\Delta u_\epsilon,\quad u_\epsilon(x,0)=u^0(x)
\end{equation}
and then put $u=\lim u_\epsilon$ as $\epsilon\to 0$. The solution of \myref{visc} should exist only
for $t\geq 0$ and so should the {\it viscosity} solution $u$. Therefore, the corresponding phase
flow is {\it irreversible} w.r.t. time.

In the classical language, this probably means  that  the set of initial conditions for a general
Coulomb system, which approach a singular set $\{x_i=x_j\}$ at finite time, has a positive Liouville
measure.

If the above picture is correct, it follows that a formally reversible Newtonian dynamics can be, in
fact, irreversible. This is a way to avoid logical contradiction between mechanics and
ther\-mo\-dy\-na\-mics at least at this particular point.
\subsection{Generalization of the Osgood conditions}
The question is: is it possible to find a proper generalization of the Osgood condition in the
spirit of the DiPerna--Lions theory. The Osgood condition (which generalizes the Lipschitz one and
guarantees the existence and uniqueness of the phase flow) is that
\begin{equation}\label{osgood}
  \int_0^1 \frac{dt}{\omega(t)}=\infty,
\end{equation}
where $\omega$ is a modulus of continuity for the vector field $b$.

For example,  $\omega(t)=t\log(1/t)$ is a typical modulus satisfying \myref{osgood}. One can show
that there exists a function $u$ with this modulus of continuity such that the distributional
derivative $\nabla u$ is not a measure. For instance, the Weierstrass function
$$u(t)=\sum_{k=1}^\infty 2^{-k}\exp(i2^{k}t)$$ is so. Indeed, the difference  $ u(t+h)-u(t)$ is
equal to
\begin{equation}\label{W}
 \sum_{n=1}^{M}2^{-n}(\exp(i2^{n}h)-1)\exp(i2^{n}t)+\sum_{n=M+1}^{\infty}2^{-n}(\exp(i2^{n}h)-1)\exp(i2^{n}t)
\end{equation}
for any $M$. We choose $M$ so that $2^Mh=o(1)$ as $h\to 0$. For instance, $$M=\log_2(1/h\log(1/h))$$
is a good choice. Now the first sum in \myref{W} can be estimated as $O(hM)=O(h\log(1/h))$, because
each term $2^{-n}(\exp(i2^{n}h)-1)$ is $O(h)$ since $2^{n}h=o(1)$, while the second sum is
$O(2^{-M})=O(h\log(1/h))$. This proves that $h\log(1/h)$ is a modulus of continuity for $u$.

If the derivative $f(t)=\sum_{k=1}^\infty \exp(i2^{k}t)$ is a measure, then by the Riesz brothers
theorem  \cite{Koosis} (if all negative Fourier coefficients of a measure vanish, then it is
absolutely continuous w.r.t. the Lebesgue measure) it is an $L^1$-function on the circle
$\R/2\pi\Z$. One can see immediately that
\begin{equation}
f(2t)=f(t)-e^{2it},\label{feq}
\end{equation}
and thus,
\begin{equation}
f(2^mt)=f(t)-\sum_{k=1}^{m} e^{2^kit}.\label{fineq}
\end{equation}
We have the equality of $L^1$-norms
$$
\int_0^{2\pi}|f(2^mt)|dt=\frac1{2^m}\int_0^{2^m2\pi}|f(t)|dt=\int_0^{2\pi}|f(t)|dt
$$
for any natural $m$. Therefore,  the $L^1$-norm of the trigonometric polynomial $\sum_{k=1}^{m}
e^{i2^kt}$  remains bounded as $m\to\infty$. This, however, contradicts the now proved Littlewood
conjecture \cite{Konyagin}, \cite{McGehee} that the $L^1$-norm of a  polynomial $\sum_{k=1}^N
a_ke^{in_kt}$ such that $|a_k|\geq 1$ and the integers $n_k$ are distinct, grows at least like $
C\log N $.
\subsection{Continuity of the phase flow}
The problem is to indicate conditions in the spirit of the DiPerna--Lions theory which guarantee the
continuity of the phase flow $X_t(x)$ w.r.t. $x$. Another face of the issue is to find {\it a
priori} Sobolev smooth\-ness for the phase flow.
\subsection{Classical interpretation of measurable phase flow} An
example of a related problem is as follows: Does it follow from the existence of measurable state
flow for the Hamiltonian \myref{ham} that classical trajectories never hit the singular set
$\{x_i=x_j\}$ for a set of initial points of full measure?

Another problem in this area is the comparison of DiPerna--Lions flows with another kind of flows
for discontinuous vector fields, like vibrosolutions, or Filippov's trajectories. \mysection{Details
and proofs}
\subsection{A priori estimates}
We note that the DiPerna--Lions theory is not totally independent of the classical theory around the
Cauchy-Lipschitz theorem. All the arguments in the DiPerna--Lions paper go via
regu\-la\-ri\-za\-tion of the Cauchy problem and then taking a limit as the small
regu\-la\-ri\-za\-tion parameter $\epsilon\to 0$. To say something about regularized problem we
utilize the classical theory. The following statement about classical solutions of \myref{cp} is
trivial.
\begin{proposition}\hglue -.3em{\bf .}\label{apriori}
Let $u^t(x)=u(x,t)$ be the classical solution of \myref{cp} at time $t$. Suppose that all data ($b$
and $u^0$) is regular. Then
\begin{equation}\label{estim}
  \Vert u^t\Vert_{L^\infty}\leq\Vert u^0\Vert_{L^\infty}.
\end{equation}
\end{proposition}
As usual, we will use this a priori estimate to construct solutions to our initial irregular Cauchy
problem, so that those solutions satisfy the same estimate. The next $L^1$ estimate of the classical
solutions is almost as easy as the previous $L^\infty$ one. We present it, in particular, for the
sake of explaining the role of the condition $\div b\in L^\infty$.
\begin{proposition}\hglue -.3em{\bf .}\label{apriori2}
Let $u^t(x)=u(x,t)$ be the classical solution of \myref{cp} at time $t\in[0,T]$. Suppose that all
data ($b$ and $u^0$) is regular. Then
\begin{equation}\label{estim2}
  \Vert u^t\Vert_{L^1}\leq C\Vert u^0\Vert_{L^1},
\end{equation}
where the constant $C$ depends only on $T$ and $M=\sup_{x\in{\mathcal T}}|\div b(x)|$.
\end{proposition}
{\it Proof}. If $u$ is a classical solution of \myref{cp} then $v=|u|$ is a weak solution. By
integrating we obtain
\begin{equation}\label{abs}
\int_{\mathcal T}v(x,t)dx-\int_{\mathcal T}v(x,0)dx=-\int_0^t\int_{\mathcal T}\div(b(x))v(x,s)dx\,ds
\end{equation}
and the modulus of the right-hand side is $\leq M \int_0^t\int_{\mathcal T}v(x,s)dx\,ds$. Now, we
have for the positive function $f(t)=\int_{\mathcal T}v(x,t)dx$ the inequality $f(t)-f(0)\leq
M\int_0^t f(s)ds$. It remains to apply the Gron\-wall lemma to get the desired estimate for
$f(t)=\Vert u^t\Vert_{L^1}$.
\subsection{Regularization}
For the sake of regularization we utilize a classical tool: convolution with a $\delta$-shaped
sequence of $C_0^\infty$ functions. More precisely, let $\rho\in C_0^\infty(\mathcal T)$ be a smooth
function with a compact support such that $\int_{\mathcal T}\rho(x)dx=1$. We assume that the support
lies within a ball on the torus, and the ball lifts homeomorphically to the universal covering
${\mathbf R}^n$. This allows to regard $\rho$ as a function on ${\mathbf R}^n$ and apply to it some
simple constructions related to ${\mathbf R}^n$. In particular, for $0<\epsilon\leq 1$ the function
$\rho_\epsilon(x)=\epsilon^{-n}\rho(x/\epsilon)$ is well defined, and  $\rho_\epsilon\to\delta_0$ in
the space of distributions as $\epsilon\to 0$. We define the convolution operator
\begin{equation}\label{conv}
  C_\epsilon u(x)=\int_{\mathcal T}u(x-y)\rho_\epsilon(y)dy,
\end{equation}
and will often write $u_\epsilon $ instead of $ C_\epsilon u$ for brevity.

The analytic heart of the DiPerna--Lions paper is a statement about commu\-ta\-tor of the operators
$C_\epsilon$ and our main differential operator (vector field)
\begin{equation}\label{diff}
  Bu=\sum_i b_i\frac{\partial u}{\partial x_i}=b\cdot\nabla u.
\end{equation}
Namely, under extended DiPerna--Lions conditions \myref{dip_l} the com\-mu\-ta\-tor $[B,C_\epsilon]$
is small in a suitable sense as $\epsilon\to 0$. This is not surprising since the operator
$C_\epsilon$ becomes arbitrary close to identity as $\epsilon\to 0$, and statements of this kind are
well known in the realm of PDE since \cite{F}.
\begin{theorem}\hglue -.3em{\bf .}\label{commut}
Suppose that the extended DiPerna--Lions conditions \myref{dip_l} hold. Then the operators
$[B,C_\epsilon]:C({\mathcal T})\to L^1({\mathcal T})$ are uniformly bounded and tend to zero
strongly. In other words, if \,$u\in C({\mathcal T})$, then the difference
$Bu_\epsilon-(Bu)_\epsilon$ has small $L^1$-norm if $\epsilon$ is small enough.
\end{theorem}
{\it Remark}. Under the original DiPerna--Lions conditions \myref{dip_le_a} the cor\-res\-pon\-ding
statement is  stronger: $[B,C_\epsilon]:L^\infty({\mathcal T})\to L^1({\mathcal T})$ is uniformly
bounded and tends to zero strongly. In other words, under original DiPerna--Lions conditions $u$ in
the proposition may not be continuous.

\noindent{\it Proof}. As usually, proof is performed in two steps: first, we prove that the
operators $[B,C_\epsilon]$ are uniformly bounded for all $\epsilon$ and $B$ subject to uniform
bounds \myref{dip_l}, second, we check that $[B,C_\epsilon]u$ is small if the vector field $b$ and
function $u$ are smooth. The second part is, in fact, trivial or, at least, well known, and we skip
it.

To prove the first part we start with the explicit formula
\begin{equation}\label{formula}
[B,C_\epsilon]u(x)=\int u(y)\{(b(y)-b(x))\cdot\nabla\rho_\epsilon(x-y)\}dy-(u\div(b))*\rho_\epsilon.
\end{equation}
We have to bound uniformly the right-hand side in $L^1$, provided that $u$ is continuous, $\Vert
u\Vert_{L^\infty}$, $\Vert \div(b)\Vert_{L^\infty}$, and $\Vert b\Vert_{W_*^{1,1}}$ are uniformly
bounded. This is trivial for the second term $(u\div(b))*~\rho_\epsilon$ because of the well known
properties of the convolution operator $C_\epsilon:L^\infty\to L^1$.

It remains to estimate
\begin{equation}\label{estimate}
  \int\left\vert\int
u(y)\{(b(y)-b(x))\cdot\nabla\rho_\epsilon(x-y)\}dy\right\vert\,dx
\end{equation}
which, as one can see after the change of variables $x=y+\epsilon z$, is not greater than
\begin{equation}\label{estimate2}
   C\int_{\mathcal B}\left\vert\int_{\mathcal T}
u(y)\frac{(b(y)-b(y+\epsilon z))}{\epsilon}dy\right\vert\,dz,
\end{equation}
where the constant $C=\sup_{z\in {\mathcal B}} |\nabla\rho(z)|$. Here, $\mathcal B$ is a small ball,
where the support of the mollifier $\rho$ is located. We can regard it as a ball in ${\mathbf R}^n$,
and the expression like $y+\epsilon z$ makes sense for $z\in {\mathcal B}$ and $y\in{\mathcal T}$.

 Now we note that for $u\in C({\mathcal T})$
\begin{equation}\label{estimate4}
  \sup_{\epsilon,z\in {\mathcal B}}\left\vert\int_{\mathcal T}u(y)\frac{(b(y)-b(y+\epsilon z))}{\epsilon}
  dx\right\vert\leq C\Vert b\Vert_{W_*^{1,1}}\Vert u\Vert_{L^\infty},
\end{equation}
where $C$ is an absolute constant, and this gives the desired estimate for the first term in
\myref{formula}.
\subsection{Approximable solutions}\label{Approximable_solutions}
\noindent{\bf Existence.} Now we apply the above regularization estimates to the Cauchy problem
\myref{cp}. We start with an ``approximate'' Cauchy problem
\begin{equation}\label{cp_delta}
\frac{\partial u_{\delta}}{\partial t}=b_\delta\cdot\nabla u_{\delta},\quad
u_{\delta}(x,0)=u^0_{\delta}(x),
\end{equation}
where the subscript $\delta$ in $b$ and $u^0$ indicates the convolution with the mollifier
$\rho_\delta$. This problem is regular and can be solved in the classical sense by using the
Cauchy-Lipschitz theorem. In particular, the functions $u_{\delta}$ are smooth. Now we consider the
``approximate solutions'' $u_{\delta,\epsilon}=u_{\delta}*\rho_\epsilon$ of the Cauchy problem
\myref{cp_delta}. Denote by $B_\delta$ the operator $b_\delta\cdot\nabla$. We have
\begin{equation}\label{cp_delta_eps}
\frac{\partial u_{\delta,\epsilon}}{\partial t}=(b_\delta\cdot\nabla
u_{\delta})*\rho_\epsilon=b_\delta\cdot\nabla u_{\delta,\epsilon}-[B_\delta,C_\epsilon]u_\delta.
\end{equation}
In view of Theorem \ref{commut} the remainder $r_{\delta,\epsilon}=-[B_\delta,C_\epsilon]u_\delta$
is uniformly w.r.t. $\delta$ small if $\epsilon$ is small. Indeed, $u_\delta$ is a continuous
function which is a priori bounded by $\Vert u^0\Vert_{L^\infty}$, while the $W_*^{1,1}$ norm of
$b_\delta$ is bounded by the $W_*^{1,1}$ norm of $b$. Therefore, we get
\begin{equation}\label{cp_delta_eps2}
\frac{\partial u_{\delta,\epsilon}}{\partial t}=b_\delta\cdot\nabla
u_{\delta,\epsilon}+r_{\delta,\epsilon},\quad u_{\delta,\epsilon}(x,0)=u^0_{\delta,\epsilon}(x),
\end{equation}
where the remainder $r_{\delta,\epsilon}\in L^\infty([0,T];L^1({\mathcal T}))$ is uniformly small.

\noindent{\bf Lipschitz bound.} Fix an $\epsilon>0$ and put $\delta\to 0$. One can see easily from
\myref{cp_delta_eps2} that the functions $[0,T]\ni t\mapsto u_{\delta,\epsilon}(t)$ are uniformly
Lipschitz as functions with values in a ``negative'' Sobolev space $H^{-s}({\mathcal T})$ for some
sufficiently large $s$. Indeed, consider the scalar product $P$ of the right-hand side of
\myref{cp_delta_eps2} with a smooth test function $\phi(x)$. We have
\begin{equation}\label{h_s}
P=\int(b_\delta\cdot\nabla u_{\delta,\epsilon}+r_{\delta,\epsilon})\phi \,dx=\int
[(-u_{\delta,\epsilon}\div b_\delta\,
+r_{\delta,\epsilon})\phi-u_{\delta,\epsilon}b_\delta\cdot\nabla \phi]\,dx
\end{equation}
In view of condition \myref{dip_l}, estimates for $r_{\delta,\epsilon}$ and $L^\infty$ a priori
estimates for $u_{\delta,\epsilon}$ one obtains the bound $|P|\leq C\Vert \phi\Vert$, where $C$ is
an absolute constant and $\Vert \phi\Vert=\sup |\phi(x)|+\sup |\nabla\phi(x)|$. Thus, the right-hand
side of \myref{cp_delta_eps2} is uniformly bounded in $H^{-s}$
 if $s$ is such that the norm $\Vert
\phi\Vert$ is continuous on $H^s$. By the Sobolev lemma we can take any $s>\frac{n}{2}+1$.

\noindent{\bf Existence (continued).} Now, we can extract a subsequence $\delta\to 0$ such that
$u_{\delta,\epsilon}(t)\to u_{\epsilon}(t)$ in $H^{-s}$ uniformly on $[0,T]$. Moreover, we can
assume that the functions $r_{\delta,\epsilon}$ converge as distributions to a measure
$r_{\epsilon}\in L^\infty([0,T];M({\mathcal T}))$ which is small with $\epsilon$. Here, $M({\mathcal
T})$ stands for the space of measures with the norm given by the total variation: $$\Vert
r\Vert_M=\sup_\phi\int r\phi,$$ where $\phi$ runs over continuous functions such that $|\phi|\leq
1$. This implies, that
\begin{equation}\label{cp_eps2}
\frac{\partial u_{\epsilon}}{\partial t}=b\cdot\nabla u_{\epsilon}+r_{\epsilon},\quad
u_{\epsilon}(x,0)=u^0_{\epsilon}(x),
\end{equation}
where both sides are in $H^{-s}$. Moreover, all the functions $u_{\epsilon}(t)$ are smooth w.r.t.
\,$x\in{\mathcal T}$, so that equation \myref{cp_eps2} is valid in the classical sense.  Now we can
again apply the same arguments on compactness in $C([0,T];H^{-s})$ and extract a subsequence
$\epsilon\to 0$ such that $u_{\epsilon}(t)\to u(t)$ in $H^{-s}$ uniformly on $[0,T]$. In other
words, we constructed an approximable solution $u$ of the Cauchy problem \myref{cp}, and the
existence part of Theorem \ref{dlth} is done.

\noindent{\bf Uniqueness.} To prove uniqueness, we have to show that an appro\-xi\-mable solution
$u$ with the zero initial condition $u^0=0$ is zero. Let $u_{\epsilon}$ be an approximate classical
solutions, satisfying \myref{cp_eps2}, where $u^0_{\epsilon}=0$. Put $v_{\epsilon}=u_{\epsilon}^2$.
One can see immediately that $v_{\epsilon}$ satisfies similar equation
\begin{equation}\label{cp_eps3}
\frac{\partial v_{\epsilon}}{\partial t}=b\cdot\nabla v_{\epsilon}+r'_{\epsilon},\quad
v_{\epsilon}(x,0)=0,
\end{equation}
where $r'_{\epsilon}=2u_{\epsilon}r_{\epsilon}$ is again a measure with a (uniformly w.r.t. \,$t$)
small with $\epsilon$ total variation. This proves, in particular, that $v=u^2$ is an approximable
solution of \myref{cp} with initial value $v^0(x)=u^0(x)^2$.

\noindent{\bf Renormalizability.} Similar arguments, where we consider
$v_{\epsilon}=\beta(u_\epsilon)$ instead of $v_{\epsilon}=u_\epsilon^2$ prove that an approximable
solution is renormalizable.

\noindent{\bf Uniqueness (continued).} Now we consider the integral $I(t)=\int v_\epsilon(x,t)dx$,
and take the integral of both sides of \myref{cp_eps3}. From the bound \myref{dip_l} on $\div b$ we
obtain immediately that $\left|\int b\cdot\nabla v_{\epsilon}(t)\right|\leq CI(t)$, where $C$ is an
absolute constant, and the integral $R_\epsilon(t)=\int_{\mathcal T} r'_{\epsilon}(t)dx$ is
uniformly small. This implies (in view of the Gronwall lemma as applied to the integral inequality
$I(t)\leq C\int_0^tI(s)ds+\int_0^tR_\epsilon(s)ds$) that $I(t)\to 0$ as $\epsilon\to 0$ uniformly
w.r.t. $t$. Therefore, $u_\epsilon\to 0$ in $L^2$, and, therefore, in $H^{-s}$. Thus, $u=0$ and the
uniqueness is done.

\noindent{\bf Regularity.} Now, the only part to be proved of Theorem \ref{dlth} is the inclusion
\myref{space}. The part $u\in L^\infty([0,T]\times{\mathcal T})$ is trivial in view of Proposition
\ref{apriori}, and it remains to show that $u\in C([0,T];L^p({\mathcal T}))$. In other words, we
have to prove that if $t_n\to t$ then $u(t_n)\to u(t)$ in $L^p$ for any $p\geq 1$.

What we already know is that $u(t_n)\to u(t)$ as distributions (even in $H^{-s}$). On the other
hand, $u(t_n)$ are uniformly bounded (in $L^\infty$). These facts combined imply that $u(t_n)\to
u(t)$ weakly in $L^2$. Since $v=u^2$ is also an approximable solution of \myref{cp} we also get that
$u(t_n)^2\to u(t)^2$ weakly in $L^2$. This allow us to prove easily that $u(t_n)\to u(t)$ (strongly)
in $L^2$. Indeed,
\begin{equation}\label{conv2}
  \int (u(t_n)-u(t))^2= \int u(t_n)^2+\int u(t)^2 -2\int
  u(t_n)u(t),
\end{equation}
where the notation $\int f=\int_{\mathcal T}f(x)dx$ is utilized. Since $u(t_n)^2\to u(t)^2$ weakly,
we conclude that $\lim \int u(t_n)^2=\int u(t)^2$, and
$$\lim\int
  u(t_n)u(t)=\int u(t)^2$$ in view of the weak convergence
  $u(t_n)\to u(t)$. Thus, $\lim\int (u(t_n)-u(t))^2=0$.
  Now, since $u$ is a uniformly bounded function, we obtain immediately that $\lim\int |u(t_n)-u(t)|^p=0$ for any $p\geq 1$, and we are done.
\subsection{Functional analysis}
Here, we prove Theorem \ref{dpth} about recovering  a measurable map from its action on bounded
measurable functions. This follows easily from the Dunford--Pettis theorem (cf. \cite{Bourb},
\cite{Schred}). We give the statement of a particular case we need.
\begin{theorem}\hglue -.3em{\bf .} (Dunford--Pettis theorem)\label{Dunford}
Suppose that
\begin{equation}\label{dp1}
  A: L^\infty({\mathcal T})\to L^\infty({\mathcal T})
\end{equation}
is a bounded operator. Then there exists a unique (modulo null sets of the Lebesgue measure) measure
$K(x,dy)$ on ${\mathcal T}$ which is measurable w.r.t. $x\in{\mathcal T}$, and such that for almost
every $x\in{\mathcal T}$
\begin{equation}\label{dp2}
  Af(x)=\int K(x,dy)f(y).
\end{equation}
\end{theorem}
We have to show that if $A$ is a homomorphism of {\it rings with unit} then the measure $K(x,dy)$ is
a $\delta$-measure: i.e., is supported by a single point $y=\Phi(x)$.

Indeed, consider the quadratic form $Q(f)=\int (Af^2(x)-Af(x)^2)dx$ which is $\equiv 0$, since $A$
is a homomorphism. In terms of the measure $K(x,dy)$ it is equal to
\begin{equation}\label{Q}
  Q(f)=\int_{{\mathcal T}}\int K(x,dy)\left[f(y)-Af(x)\right]^2dx\equiv 0.
\end{equation}
Therefore, $$\int K(x,dy)\left[f(y)-Af(x)\right]^2=0$$ for almost all $x$. For any such an
  $x$ we get that any function $f$ takes a single value $Af(x)=\int
  K(x,dy)f(y)$ modulo null sets for the measure $K(x,dy)$. This
  means that the support of $K(x,dy)$ consists of a single
  point which is exactly what we have to prove.
\subsection{Stability theorem}
Theorem \ref{stabth} can be proved by methods already used in the proof of Theorem \ref{dlth}.

We will show that there is a convergent subsequence $u_n\to u$ in $C([0,T];L^1)$ to an approximable
solution $u$. In view of the uniqueness of the approximable solution this will prove the Theorem.

Let $u_{n,\epsilon}$ be an approximate solution to the Cauchy problem \myref{cp} with $b$, resp.
$u^0$ replaced by $b_n$ resp. $u^0_n$. More precisely, we assume that
\begin{equation}\label{appr}
\frac{\partial u_{n,\epsilon}}{\partial t}=b_n\cdot\nabla u_{n,\epsilon}+r_{n,\epsilon},\quad
u_{n,\epsilon}(x,0)=u^0_{n,\epsilon}(x),
\end{equation}
where $\Vert r_{n,\epsilon}\Vert_M\leq\epsilon$, and $\Vert
u^0_{n}-u^0_{n,\epsilon}\Vert_{L^1}\leq\epsilon$. The arguments utilized in the proof of the {\bf
Lipshitz bound} in the previous subsection \ref{Approximable_solutions} show that the right-hand
side of \myref{appr} is uniformly bounded in $H^{-s}$ for any $s>\frac{n}{2}+1$. This, in turn,
shows that for a subsequence of indices $n$ we have $u_{n,\epsilon}\to u_{\epsilon}$ in
$C([0,T];H^{-s})$ and the functions $u_{\epsilon}$ satisfy
\begin{equation}\label{cp_eps4}
\frac{\partial u_{\epsilon}}{\partial t}=b\cdot\nabla u_{\epsilon}+r_{\epsilon},\quad
u_{\epsilon}(x,0)=u^0_{\epsilon}(x),
\end{equation}
where $\Vert r_{\epsilon}\Vert_M\leq\epsilon$, $\Vert u^0-u^0_{\epsilon}\Vert_{L^1}\leq\epsilon$.

Now we can again apply the same arguments on compactness in $C([0,T];H^{-s})$ and extract a
subsequence $\epsilon\to 0$ such that $u_{\epsilon}(t)\to u(t)$ in $H^{-s}$ uniformly on $[0,T]$.

This function $u$ is an approximable solution to \myref{cp}, and $u_\epsilon$ are its
approximations.

Moreover, for a subsequence of indices $n$ we have $u_{n}(t)\to u(t)$ in $H^{-s}$ uniformly on
$[0,T]$.

Now, an easy adaptation of the arguments about {\bf Regularity} from subsection
\ref{Approximable_solutions} shows that $u_{n}(t)\to u(t)$ in $L^p$ for any $p<\infty$, $t\in[0,T]$,
and $u\in C([0,T];L^p)$ for any $p<\infty$. This implies $u_{n}\to u$ in $C([0,T];L^p)$ and we are
done.
\def\refname{References}

\end{document}